\documentclass{amsart}
\usepackage{amssymb}

\newtheorem{thm}{Theorem}
\newtheorem{cor}{Corollary}
\newtheorem{lem}{Lemma}

\newtheorem{rem}{Remark}

\newtheorem{conj}{Conjecture}
\theoremstyle{definition}
\newtheorem{example}[equation]{Example}
\newtheorem{prob}[equation]{Problem}

\newcommand{\D}{{\mathbb D}}

\def\be{\begin{equation}}
\def\ee{\end{equation}}

\newcommand{\bee}{\begin{enumerate}}
\newcommand{\eee}{\end{enumerate}}

\newcommand{\blem}{\begin{lem}}
\newcommand{\elem}{\end{lem}}
\newcommand{\bthm}{\begin{thm}}
\newcommand{\ethm}{\end{thm}}
\newcommand{\bcor}{\begin{cor}}
\newcommand{\ecor}{\end{cor}}
\newcommand{\beg}{\begin{example}}
\newcommand{\eeg}{\end{example}}
\newcommand{\begs}{\begin{examples}}
\newcommand{\eegs}{\end{examples}}
\newcommand{\bdefe}{\begin{defin}}
\newcommand{\edefe}{\end{defin}}
\newcommand{\bprob}{\begin{prob}}
\newcommand{\eprob}{\end{prob}}
\newcommand{\bei}{\begin{itemize}}
\newcommand{\eei}{\end{itemize}}

\newcommand{\bcon}{\begin{conj}}
\newcommand{\econ}{\end{conj}}
\newcommand{\bcons}{\begin{conjs}}
\newcommand{\econs}{\end{conjs}}
\newcommand{\bprop}{\begin{propo}}
\newcommand{\eprop}{\end{propo}}
\newcommand{\br}{\begin{rem}}
\newcommand{\er}{\end{rem}}
\newcommand{\brs}{\begin{rems}}
\newcommand{\ers}{\end{rems}}
\newcommand{\bo}{\begin{obser}}
\newcommand{\eo}{\end{obser}}
\newcommand{\bos}{\begin{obsers}}
\newcommand{\eos}{\end{obsers}}
\newcommand{\bpf}{\begin{pf}}
\newcommand{\epf}{\end{pf}}
\newcommand{\ba}{\begin{array}}
\newcommand{\ea}{\end{array}}
\newcommand{\beq}{\begin{eqnarray}}
\newcommand{\beqq}{\begin{eqnarray*}}
\newcommand{\eeq}{\end{eqnarray}}
\newcommand{\eeqq}{\end{eqnarray*}}

\begin{document}
\bibliographystyle{amsplain}

\title[On the difference of initial logarithmic coefficients ]{On the difference of initial logarithmic coefficients for the class of univalent functions}

\author[M. Obradovi\'{c}]{Milutin Obradovi\'{c}}
\address{Department of Mathematics,
Faculty of Civil Engineering, University of Belgrade,
Bulevar Kralja Aleksandra 73, 11000, Belgrade, Serbia}
\email{obrad@grf.bg.ac.rs}

\author[N. Tuneski]{Nikola Tuneski}
\address{Department of Mathematics and Informatics, Faculty of Mechanical Engineering, Ss. Cyril and
Methodius
University in Skopje, Karpo\v{s} II b.b., 1000 Skopje, Republic of North Macedonia.}
\email{nikola.tuneski@mf.edu.mk}

\subjclass{30C45, 30C50, 30C55}
\keywords{univalent, Grunsky coefficients,  logarithmic coefficients,  difference}

\begin{abstract}
In this paper we give estimates of the differences $|\gamma _{3}|-|\gamma_{2}|$ and $|\gamma _{4}|-|\gamma_{3}|$ for the class of functions $f$ univalent in the unit disc and normalized by $f(0)=f'(0)-1=0$. Here, $\gamma_{2}$, $\gamma_{3}$ and $\gamma_{4}$ are the initial logarithmic coefficients of the function $f$.
\end{abstract}

\maketitle

\medskip

\section{Introduction and preliminaries}
As usual, let $\mathcal{A}$ be the class of functions $f$ that are analytic  in the open unit disc $\D=\{z:|z|<1\}$ of the form
\begin{equation}\label{eq-1}
f(z)=z+a_{2}z^{2}+a_{3}z^{3}+\cdots,
\end{equation}
and let $\mathcal{S}$ be the subclass of $\mathcal{A}$ consisting of functions that are univalent in $\D$.

One of the problems in the theory of univalent functions is finding sharp estimates of logarithmic coefficient, $\gamma_n$, of a univalent function $f \in\mathcal{S}$,  defined by
\begin{equation}\label{eq-2}
   \log\frac{f(z)}{z}=2\sum_{n=1}^\infty \gamma_{n} z^{n}.
\end{equation}
From the relations \eqref{eq-1} and \eqref{eq-2}, after  equating the coefficients we receive the next initial logarithmic coefficients:
\begin{equation}\label{eq-3}
 \gamma_{1}=\frac{a_{2}}{2},\quad \gamma_{2}=\frac{1}{2}\left(a_3-\frac{1}{2}a_{2}^{2}\right), \quad \text{and}\quad \gamma_3 =\frac{1}{2}\left(a_4-a_2a_3+\frac13a_2^3\right).
\end{equation}
Relatively little exact information is known about the logarithmic coefficients. The natural conjecture $|\gamma_n|\le1/n$, inspired by the Koebe function (whose logarithmic coefficients are $1/n$) is false even in order of magnitude (see Duren \cite[Section 8.1]{duren}), and true only for the class of starlike functions (\cite{der}). For the class $\mathcal{S}$ the sharp estimates of single logarithmic coefficients are known only for $\gamma_1$ and $\gamma_2$, namely,
$$|\gamma_1|\le1\quad\mbox{and}\quad |\gamma_2|\le \frac12+\frac1e=0.635\ldots.$$
In their papers \cite{OT_2021-3} and \cite{OT_2022-3}, the authors gave estimates $|\gamma_3|\leq 0.5566178\ldots$ and
$|\gamma_4|\leq0.51059\ldots$ for the  class $\mathcal{S}.$

\medskip

In this paper we will give the estimates for the differences $|\gamma _{3}|-|\gamma_{2}|$ and $|\gamma _{4}|-|\gamma_{3}|$ in the class $\mathcal{S}$.  The sharp estimates of $|\gamma _{2}|-|\gamma_{1}|$ were established in \cite{lecko}, with a simple proof given in \cite{OT_2023-3} using different technique.

In this paper our main tool will be a method based on Grunsky's coefficients, used previously in \cite{OT_2021-3} and \cite{OT_2022-3}. This method is different than the method used in\cite{lecko}, or in \cite{OT_2023-3}. 

We will use mainly the notations and results given in the book of N.A. Lebedev (\cite{Lebedev}).

Let $f \in \mathcal{S}$ and let
\[ \log\frac{f(t)-f(z)}{t-z}=\sum_{p,q=0}^{\infty}\omega_{p,q}t^{p}z^{q}, \]
where $\omega_{p,q}$ are so-called Grunsky's coefficients with property $\omega_{p,q}=\omega_{q,p}$.
For those coefficients we have the next Grunsky's inequality (\cite{duren,Lebedev}):
\begin{equation}\label{eq 4}
\sum_{q=1}^{\infty}q \left|\sum_{p=1}^{\infty}\omega_{p,q}x_{p}\right|^{2}\leq \sum_{p=1}^{\infty}\frac{|x_{p}|^{2}}{p},
\end{equation}
where $x_{p}$ are arbitrary complex numbers such that last series converges.

Further, it is well-known that if function
\begin{equation}\label{eq 5}
f(z)=z+a_{2}z^{2}+a_{3}z^{3}+\cdots
\end{equation}
belongs to $\mathcal{S}$, then also
\[
f_{2}(z)=\sqrt{f(z^{2})}=z +c_{3}z^3 +c_{5}z^{5}+\cdots
\]
belongs to the class $\mathcal{S}$. For the function $f_{2}$ we have the appropriate Grunsky's
coefficients of the form $\omega_{2p-1,2q-1}$, and the inequality \eqref{eq 4} has the form
\begin{equation}\label{eq-6}
\sum_{q=1}^{\infty}(2q-1) \left|\sum_{p=1}^{\infty}\omega_{2p-1,2q-1}x_{2p-1}\right|^{2}\leq \sum_{p=1}^{\infty}\frac{|x_{2p-1}|^{2}}{2p-1}.
\end{equation}
As it has been shown in \cite[p.57]{Lebedev}, if $f$ is given by \eqref{eq-1} then the coefficients $a_{2}$, $a_{3}$, $a_{4}$ and $a_5$
are expressed by Grunsky's coefficients  $\omega_{2p-1,2q-1}$ of the function $f_{2}$ given by
\eqref{eq 5} in the following way:
\begin{equation}\label{eq 7}
\begin{split}
a_{2}&=2\omega _{11},\\
a_{3}&=2\omega_{13}+3\omega_{11}^{2}, \\
a_{4}&=2\omega_{33}+8\omega_{11}\omega_{13}+\frac{10}{3}\omega_{11}^{3},\\
a_{5}&=2\omega_{35}+8\omega_{11}\omega_{33}+5\omega_{13}^{2}+18\omega_{11}^{2}\omega_{13}+\frac{7}{3}\omega_{11}^{4},\\
0&=3\omega_{15}-3\omega_{11}\omega_{13}+\omega_{11}^{3}-3\omega_{33},\\
0&=\omega_{17}-\omega_{35}-\omega_{11}\omega_{33}-\omega_{13}^{2}+\frac{1}{3}\omega_{11}^{4}.
\end{split}
\end{equation}
We note that in the book of Lebedev \cite{Lebedev} there exists a typing mistake for the coefficient $a_{5}$. Namely, instead of the term $5\omega_{13}^{2}$, there is $5\omega_{15}^{2}$.

\medskip

\section{Maine results}
\bthm\label{th-1}
Let  $\gamma_{2}$, $\gamma_{3}$ and $\gamma_{4}$ be the logarithmic coefficients of function $f\in \mathcal{S}$. Then
$$|\gamma _{3}|-|\gamma_{2}|\leq \frac{1}{\sqrt{5}} \quad\text{and}\quad
|\gamma _{4}|-|\gamma_{3}|\leq \frac{1}{\sqrt{7}}.$$
\ethm
\begin{proof} For our consideration we need the following facts.
From \eqref{eq-6}, choosing $x_{2p-1}=0$ when $p=3,4,\ldots$, we have
\[
\begin{split}
&|\omega_{11}x_{1}+\omega_{31}x_{3}|^{2}+3|\omega_{13}x_{1}+\omega_{33}x_{3}|^{2}+
5|\omega_{15}x_{1}+\omega_{35}x_{3}|^{2}\\
+&7|\omega_{17} x_1 +\omega_{37} x_3 |^2\leq |x_1|^2+\frac{|x_3|^2}{3}.
\end{split}
\]
If additionally, $x_{1}=1$ and $x_{3}=0$, then
$$|\omega_{11}|^{2}+3|\omega_{13}|^{2}+5|\omega_{15}|^{2}+7|\omega_{17}|^{2}\leq 1,$$
and from here also
$$|\omega_{11}|^{2}+3|\omega_{13}|^{2}+5|\omega_{15}|^{2}\leq 1,$$
$$|\omega_{11}|^{2}+3|\omega_{13}|^{2}\leq 1,$$
and $$|\omega_{11}|^{2}\leq1 .$$
From the previous inequalities we have
\begin{equation}\label{eq-9}
\begin{split}
|\omega_{11}|&\leq1,\\
|\omega_{13}|&\leq\frac{1}{\sqrt{3}}\sqrt{1-|\omega_{11}|^{2}}, \\
|\omega_{15}|&\leq \frac{1}{\sqrt{5}}\sqrt{1-|\omega_{11}|^{2}-3|\omega_{13}|^{2}},\\
|\omega_{17}|&\leq \frac{1}{\sqrt{7}}\sqrt{1-|\omega_{11}|^{2}-3|\omega_{13}|^{2}-5|\omega_{15}|^{2}}.
\end{split}
\end{equation}
Also, from the fifth relation in \eqref{eq 7} we obtain
\begin{equation}\label{eq-11}
\omega_{33}=\omega_{15}-\omega_{11}\omega_{13}+\frac{1}{3}\omega_{11}^{3}.
\end{equation}

\medskip

Now, for the first estimate of the theorem, using \eqref{eq-3} and \eqref{eq 7}, we have
\begin{equation}\label{eq-12}
\begin{split}
|\gamma _{3}|-|\gamma_{2}|&=\frac{1}{2}\left|a_4-a_2a_3+\frac13a_2^3\right|-
\frac{1}{2}\left|a_3-\frac{1}{2}a_{2}^{2}\right|\\
&=|\omega_{33}+2\omega_{11}\omega_{13}|-\left|\omega_{13}+\frac{1}{2}\omega_{11}^{2}\right|,
\end{split}
\end{equation}
and  after applying $|\omega_{11}|\leq1$, brings
\[
\begin{split}
|\gamma _{3}|-|\gamma_{2}|& \leq|\omega_{33}+2\omega_{11}\omega_{13}|-|\omega_{11}|\cdot \left|\omega_{13}+\frac{1}{2}\omega_{11}^{2}\right|\\
&\leq \left|(\omega_{33}+2\omega_{11}\omega_{13})-(\omega_{11})\left(\omega_{13}+\frac{1}{2}\omega_{11}^{2}\right)\right|\\
&=\left|\omega_{33}+\omega_{11}\omega_{13}-\frac{1}{2}\omega_{11}^{3}\right|.
\end{split}
\]
From here, having in mind the relation \eqref{eq-11} and inequalities \eqref{eq-9}, we receive
\begin{equation}\label{eq-14}
\begin{split}
|\gamma _{3}|-|\gamma_{2}|&\leq \left|\omega_{15}-\frac{1}{6}\omega_{11}^{3}\right| \leq |\omega_{15}|+\frac{1}{6}|\omega_{11}|^{3}\\
&\leq \frac{1}{\sqrt{5}}\sqrt{1-|\omega_{11}|^{2}-3|\omega_{13}|^{2}}+\frac{1}{6}|\omega_{11}|^{3}\\
&=:\Phi(|\omega_{11}|,|\omega_{13}|),
\end{split}
\end{equation}
where
\[\Phi(u,v)=\frac{1}{\sqrt{5}}\sqrt{1-u^{2}-3v^{2}}+\frac{1}{6}u^{3}\]
with domain
\[ \Omega = \left\{ (u,v) : 0\leq u \leq1,\, 0\leq v \leq \frac{1}{\sqrt{3}}\sqrt{1-u^{2}} \right\}.\]

\medskip

It remains to find the maximal value of $\Phi$ on the domain $\Omega$.

\medskip

Is easy to verify that  the function $\Phi(u,v)$  has no interior singular points in the domain $\Omega$ ($\Phi'_v(u,v)=0$, if, and only of, $v=0$).
On the edges:
\[
\begin{split}
\Phi(u,0)&=\frac{1}{\sqrt{5}}\sqrt{1-u^{2}}+\frac{1}{6}u^{3}\leq \Phi(0,0)= \frac{1}{\sqrt5}=0.44721\ldots ,\\
\Phi(0,v)&=\frac{1}{\sqrt{5}}\sqrt{1-3v^{2}}\leq\frac{1}{\sqrt{5}}=0.44721\ldots\\
\Phi\left(u,\frac{1}{\sqrt{3}}\sqrt{1-u^{2}}\right)&=\frac{1}{6}u^3 \leq \frac{1}{6}=0.1(6).
\end{split}
\]

\medskip

Using all previous facts and \eqref{eq-14}, we finally conclude that
$$|\gamma _{3}|-|\gamma_{2}|\leq \frac{1}{\sqrt{5}}.$$

\medskip

As for second estimate of this theorem, from \eqref{eq-2} we receive
\[\gamma_4 =\frac{1}{2}\left(a_5-a_2a_4-\frac{1}{2}a_3^2+a_{2}^{2}a_{3}-\frac{1}{4}a_{2}^{4}\right),\]
and by using the relations \eqref{eq 7},
\begin{equation}\label{eq-16}
\gamma_4 =\frac{1}{2}\left(2\omega_{35}+3\omega_{13}^{2}+4\omega_{11}\omega_{33}+4\omega_{11}^{2}\omega_{13} -\frac{5}{6}\omega_{11}^{4}\right).
\end{equation}
Next, from the last two relations in \eqref{eq 7} we can express $\omega_{33}$ and $\omega_{35}$ and apply them in
\eqref{eq-16}.  After some calculations, we receive
\[\gamma_4=\omega_{17}+\omega_{11}\omega_{15}+\omega_{11}^{2}\omega_{13}+\frac{1}{2}\omega_{13}^{2}+\frac{1}{4}\omega_{11}^{4}.\]
Now, using the above expression for $\gamma_4$,  the expression for $\gamma_{3}$  obtained within  \eqref{eq-12}, $\gamma_3=\omega_{33} + 2\omega_{11}\omega_{13}$, and the expressions for $\omega_{33}$ and $\omega_{35}$ from the last two lines from \eqref{eq 7}, having in mind that $|\omega_{11}|\leq1$, we receive:
\begin{equation}\label{eq-18}
\begin{split}
|\gamma _{4}|-|\gamma_{3}|&\leq |\gamma _{4}|-|\omega_{11}||\gamma_{3}|\\
&\leq|\gamma _{4}-\omega_{11}\gamma_{3}|\\
&=\left|\omega_{17}+\frac{1}{2}\omega_{13}^{2}-\frac{1}{12}\omega_{11}^{4}\right|\\
&\leq |\omega_{17}|+\frac{1}{2}|\omega_{13}|^{2}+\frac{1}{12}|\omega_{11}|^{4}\\
&\leq \frac{1}{\sqrt{7}}\sqrt{1-|\omega_{11}|^{2}-3|\omega_{13}|^{2}-5|\omega_{15}|^{2}}
+\frac{1}{2}|\omega_{13}|^{2}+\frac{1}{12}|\omega_{11}|^{4}\\
&\leq \frac{1}{\sqrt{7}}\sqrt{1-|\omega_{11}|^{2}-3|\omega_{13}|^{2}}
+\frac{1}{2}|\omega_{13}|^{2}+\frac{1}{12}|\omega_{11}|^{4}\\
&=: \Psi (|\omega_{11}|^{2},|\omega_{13}|^{2}),
\end{split}
\end{equation}
where
\[\Psi (s,t)= \frac{1}{\sqrt{7}}\sqrt{1-s-3t}+\frac{1}{2}t+\frac{1}{12}s^{2},\]
$0\leq s \leq1$, $0\leq t \leq \frac{1}{3}(1-s)$.
It is easily verified that the function $\Psi (s,t)$ attains its maximum  $\frac{1}{\sqrt{7}}$ for $s=t=0$.
Finally, from \eqref{eq-18},
$$|\gamma _{4}|-|\gamma_{3}|\leq \frac{1}{\sqrt{7}}.$$
\end{proof}

Previous theorem leads to the conjecture for the difference of the moduli of two consecutive logarithmic coefficients of univalent functions.

\begin{conj}
If $\gamma_{n}$ is logarithmic coefficient of function $f\in \mathcal{S}$, then
$$|\gamma _{n}|-|\gamma_{n-1}|\leq \frac{1}{\sqrt{2n-1}},\quad n=3,4,\ldots .$$
\end{conj}

\end{document}